\documentclass[journal]{IEEEtran}

\usepackage{algorithm}
\usepackage{algorithmic}
\usepackage{amssymb}
\usepackage{amsfonts}
\usepackage{amsmath}
\usepackage{amsthm}
\usepackage{bm}
\usepackage{booktabs}
\usepackage{cases}
\usepackage{cite}
\usepackage{color}
\usepackage{comment}
\usepackage{empheq}
\usepackage{enumerate}
\usepackage{enumitem}
\usepackage{framed}
\usepackage[pdftex]{graphicx}
\usepackage[hidelinks]{hyperref}
\usepackage{latexsym}
\usepackage{matlab-prettifier}
\usepackage{tikz}
\usepackage{url}
\usepackage{blkarray}

\newtheorem{theorem}{Theorem}
\newtheorem*{assumption}{Assumption}
\newtheorem{definition}{Definition}
\newtheorem{remark}{Remark}
\newtheorem{lemma}{Lemma}
\newtheorem{corollary}{Corollary}

\newtheorem{prop}{Proposition}

\newtheorem{problem}{Problem}

\newcommand{\tblack}{\textcolor{black}}

\newcommand{\minimize}{\mathop{\rm minimize}\limits}
\newcommand{\subjectto}{\mathop{\rm subject\ to}\limits}

\newcommand{\paren}[1]{\left(#1\right)}
\newcommand{\sbra}[1]{\left[#1\right]}

\DeclareMathOperator{\rank}{rank}

\DeclareMathOperator{\tr}{tr}

\newcommand{\Real}{\mathbb{R}}

\begin{document}
\title{Data-driven $h^{2}$ model reduction for linear discrete-time systems}
\author{Hiroki Sakamoto and Kazuhiro Sato\thanks{H. Sakamoto and K. Sato are with the Department of Mathematical Informatics, Graduate School of Information Science and Technology, The University of Tokyo, Tokyo 113-8656, Japan, email: soccer-books0329@g.ecc.u-tokyo.ac.jp (H. Sakamoto), kazuhiro@mist.i.u-tokyo.ac.jp (K. Sato) }}

\maketitle
\thispagestyle{empty}
\pagestyle{empty}

\begin{abstract}
We present a data-driven framework for $h^{2}$-optimal model reduction for linear discrete-time systems. Our main contribution is to create optimal reduced-order models in the $h^{2}$-norm sense directly from the measurement data alone, without using the information about the original system. In particular, we focus on the fact that the gradients of the $h^{2}$ model reduction problem are expressed using the discrete-time Lyapunov equation and the discrete-time Sylvester equation, and derive the data-driven gradients. 
\tblack{The proposed algorithm uses the output of an existing MOR as the initial point, and convergence to a stationary point is guaranteed under certain assumptions.}
In numerical experiments, we demonstrate that, for a modeling task in neuroscience, our method constructs a reduced-order model that outperforms DMDc in terms of the $h^2$-norm. 
\end{abstract}


\IEEEpeerreviewmaketitle

\section{Introduction} \label{sec:intro}
\tblack{Model Order Reduction (MOR) is a powerful tool for compressing high-dimensional dynamical systems, enabling faster simulation and control without sacrificing essential dynamics.}
In systems and control theory, MOR generates reduced-order models (ROMs) that are capable of reproducing the input-output behavior of large-scale dynamical systems with high accuracy.
The resulting ROMs have the advantage of being easier to manipulate and control, unlike the inherently larger systems.
There are MOR methods based on Singular Value Decomposition (SVD) \cite{moore1981principal, mullis1976synthesis, gugercin2004survey, antoulas2005approximation, willcox2002balanced}, Krylov subspaces (or moment matching methods) \cite{bai2002krylov, antoulas2005approximation, gugercin2008h_2, antoulas2010interpolatory}, and the optimization of the $h^{2}$ (or $H^{2}$) norm \cite{gugercin2008h_2, van2008h2, bunse2010h2, sato2018riemannian, sato2019riemannian}.
These classical model reduction methods based on the state-space description of the system are known as model-based MOR.
On the other hand, when a state-space description of the system is not available or computational modeling is difficult, it is desirable to create ROMs using the measurement data alone.

Data-driven model reduction constructs the ROMs that are capable of reproducing the input-output behavior of large-scale dynamical systems with high accuracy directly from the measurement data alone.
Several well-known data-driven MOR methods have been proposed due to the increasing use of data.
\tblack{The data-driven balanced truncation (BT)~\cite{rapisarda2011identification, ma2011reduced, gosea2022data, burohman2023data} is one of the common data-driven model reduction inspired by the BT, which is a model-based MOR. 
In this method, BT is performed by using the measurement data. 
For example, \cite{gosea2022data} performs MOR by estimating the Gramians from the data information.
Another well-known approach is data-driven interpolatory model reduction. This approach includes the Loewner framework
~\cite{mayo2007framework, gosea2021data} 
and the AAA algorithm~\cite{nakatsukasa2018aaa,gosea2021algorithms}.
Other data-driven methods, such as 
dynamic mode decomposition (DMD) \cite{kutz2016dynamic,proctor2016dynamic} or operator inference \cite{peherstorfer_data-driven_2016, kramer2024learning} can also be employed.}
These approaches are the leading data-driven model reduction, but they may not yield optimal ROMs in the $h^{2}$-norm sense. In such cases, the output of the resulting reduced system may not approximate that of the original system. 

In this paper, we propose a data-driven model reduction method for linear discrete-time systems. 
We focus on obtaining the optimal ROMs in the $h^{2}$-norm sense, using the measurement time-domain data. The contributions of this research are:

\noindent 1) 
We propose a new framework for data-driven model reduction for linear discrete-time systems. We focus on obtaining the ROMs that perform well in terms of the $h^2$-norm, directly from the measurement data. 
Furthermore, under certain assumptions about the measurement data, we derive the gradients of the optimization problem characterized by the data.

\noindent 2) 
We develop an optimization algorithm for the proposed method.
The proposed algorithm can be hybridized with existing MOR.
\tblack{Furthermore, convergence to a stationary point is guaranteed under certain assumptions.}
We demonstrate that, for a modeling task in neuroscience, our method constructs a ROM that outperforms DMDc in terms of the $h^2$-norm. 

The remainder of the paper is structured as follows. Section \ref{sec:pre} describes the model-based $h^{2}$ MOR problem for the discrete-time systems. We describe the problem setting of this study in Section \ref{sec:formulation}. Section \ref{sec:main} describes the data-driven $h^{2}$ MOR method proposed in this paper and its algorithm is explained in Section \ref{sec:algo}. Section \ref{subsec:computational_neuroscience} presents the results of numerical experiments.
We conclude and discuss future work in Section \ref{sec:conclusion}.

{\it Notation}:
We denote the imaginary unit by $\mathrm{i}$. 
We also denote the 2-norm of a vector $a$, and for a matrix $A$ its Frobenius norm, transpose, complex conjugate transpose (for complex $A$), trace, Moore--Penrose pseudoinverse, $i$th eigenvalue, $i$th row, and $(i,j)$th entry
by $\|a\|$, $\|A\|_{F}$, $A^{\top}$, $A^{*}$, $\tr A$, $A^{\dagger}$, $\lambda_{i}(A)$, $A(i,:)$, and $[A]_{i,j}$, respectively.

\section{Preliminaries} \label{sec:pre}
In this section, we summarize a model-based $h^{2}$ MOR \cite{antoulas2005approximation, gugercin2008h_2,van2008h2,bunse2010h2} for linear time-invariant discrete-time systems
\begin{empheq}[left=\empheqlbrace]{align} \label{eq:lti}
    \begin{aligned}
    x_{k+1} &= Ax_{k} + Bu_{k}, \\
    y_{k} &= Cx_{k}
    \end{aligned}
\end{empheq}
with transfer function
$H(z) := C(zI_{n}-A)^{-1}B$,
\tblack{where $z\in \mathbb{C}$ is the $z$-transform variable, }
$k$ is the time index, and $u_{k}\in \Real^{m}$, $x_{k}\in \Real^{n}$, and $y_{k}\in \Real^{p}$ denote the input, state, and output at time $t_{k}$, respectively.
The matrices $A\in \Real^{n\times n}$, $B\in \Real^{n\times m}$, and $C\in \Real^{p\times n}$ are constant matrices, which are known.
Furthermore, the $h^{2}$-norm of the system~\eqref{eq:lti} is defined as
\begin{align*}
    \|H\|_{h^{2}}:= \paren{\frac{1}{2\pi} \int_{0}^{2\pi}\tr \sbra{H^{*}(e^{-\mathrm{i}\theta})H(e^{\mathrm{i}\theta})}d\theta}^{\frac{1}{2}}.
\end{align*}

In this study, we assume that $C = I_{n}$. This means $p=n$. That is, we consider a system in which all state data are observed as outputs. This assumption is commonly made in the context of DMD, as explained in Section \ref{sec:formulation}.
In addition, we assume that the system~\eqref{eq:lti} is asymptotically stable, i.e. all eigenvalues of $A\in \Real^{n\times n}$ lie inside the unit circle. 
Such a matrix $A$ is called a stable matrix.

The ROM of the system~\eqref{eq:lti} is defined as
\begin{empheq}[left=\empheqlbrace]{align}\label{eq:reduced_lti}
    \begin{aligned}
    \hat{x}_{k+1} &= \hat{A}\hat{x}_{k} + \hat{B}u_{k}, \\
    \hat{y}_{k} &= \hat{C}\hat{x}_{k}
    \end{aligned}
\end{empheq}
with transfer function $\hat{H}(z) := \hat{C}(zI_{r}-\hat{A})^{-1}\hat{B}$,
where $\hat{x}_{k}\in \Real^r$, $\hat{y}_{k}\in \Real^{n}$, and $r\ll n$.
\tblack{
Following \cite{gugercin2008h_2,sato2019riemannian}, we obtain}
\begin{equation}\label{eq:relation_y_h2_u}
  \tblack{\max_{k\in \mathbb{Z}} \|x_k-\hat{C}\hat{x}_k\|
  \leq
  \|H-\hat H\|_{h^{2}}\cdot \sqrt{\sum_{k\in \mathbb{Z}} \|u_k\|^2}.}
\end{equation}
\tblack{\eqref{eq:relation_y_h2_u} means that, under sufficiently small $\|u_k\|$, the output error $\|x_k - \hat{C}\hat{x}_k\|$ becomes smaller by minimizing the approximation error $\|H-\hat H\|_{h^{2}}$.
}

In this context, $\hat{C}$ can be regarded as a matrix that reconstructs the dimensionally reduced state data back into the original output space, and it is employed to compare the original system~\eqref{eq:lti} and ROM~\eqref{eq:reduced_lti}.

\subsection{Model-Based \texorpdfstring{$h^2$}{h2} Model Order Reduction}
Motivated by \eqref{eq:relation_y_h2_u}, we focus on finding a ROM \eqref{eq:reduced_lti}
that minimizes $\|H-\hat H\|_{h^{2}}$ so as to bound the output error $\|x_k-\hat{C}\hat{x}_k\|$,
under the constraint $(\hat{A},\hat{B},\hat{C})\in\mathcal{E}_s$, where
\begin{align}\label{eq:feasible_set}
    \tblack{\mathcal{E}_s:=\{(\hat{A},\hat{B},\hat{C})\in \mathbb{R}^{r\times r}\times \mathbb{R}^{r\times m}\times \mathbb{R}^{n\times r}\mid \hat{A}\ \text{is stable}\}.}
\end{align}

\begin{problem}\label{prob:h2}
The $h^{2}$ MOR problem for \eqref{eq:lti} is
\begin{equation*}
\begin{aligned}
    &\minimize && \|H-\hat{H}\|_{h^{2}} \\
    &\subjectto && (\hat{A}, \hat{B}, \hat{C}) \in \mathcal{E}_s.
\end{aligned}
\end{equation*}
\end{problem}

\begin{prop}[\!\cite{antoulas2005approximation}\cite{bunse2010h2}]
\label{prop:modified_opt}
    Consider the asymptotically stable systems \eqref{eq:lti} with $C=I_n$ and \eqref{eq:reduced_lti}.
    Then, the objective function of Problem~\ref{prob:h2} can be rewritten as
    \begin{align*}
        \|H-\hat{H}\|_{h^{2}}^{2} &= \tr (\Sigma_{c}) + f(\hat{A}, \hat{B}, \hat{C}) \\
        &= \tr (B^{\top}\Sigma_{o}B) + f(\hat{A}, \hat{B}, \hat{C}).
    \end{align*}
    Here, 
    \begin{align*}
        f(\hat{A},\hat{B}, \hat{C})&:= \tr(\hat{C}P\hat{C}^{\top}) - 2\tr (R\hat{C}^{\top}) \\
        &=\tr (\hat{B}^{\top}Q\hat{B}) + 2\tr (B^{\top}S\hat{B}),
    \end{align*}
    and $\Sigma_{c}$, $\Sigma_{o}$, $P$, $Q$, $R$, and $S$ are the unique solutions of the following equations:
    \begin{align}
        A\Sigma_{c}A^{\top} + BB^{\top} &= \Sigma_{c}, \label{eq:LYAP1}\\
        A^{\top}\Sigma_{o}A + I_{n} &= \Sigma_{o}, \label{eq:LYAP2}\\
        \hat{A}P\hat{A}^{\top} + \hat{B}\hat{B}^{\top} &= P, \label{eq:lyap1}\\
        \hat{A}^{\top}Q\hat{A} + \hat{C}^{\top}\hat{C} &= Q, \label{eq:lyap2}\\
        AR\hat{A}^{\top} + B\hat{B}^{\top} &= R, \label{eq:sylve1}\\
        A^{\top}S\hat{A} - \hat{C} &= S, \label{eq:sylve2}
    \end{align}
    where \eqref{eq:lyap1} and \eqref{eq:lyap2} are the discrete-time Lyapunov equations and \eqref{eq:sylve1} and \eqref{eq:sylve2} are the discrete-time Sylvester equations.
\end{prop}
From Proposition~\ref{prop:modified_opt}, Problem~\ref{prob:h2} can be rewritten as
\begin{align}
\label{eq:opt_modified}
\begin{aligned}
    &\minimize &&f(\hat{A}, \hat{B}, \hat{C}) \\
    &\subjectto && (\hat{A}, \hat{B}, \hat{C}) \in \mathcal{E}_s.
\end{aligned}
\end{align}
A modification as in \eqref{eq:opt_modified} eliminates the need to compute the large Lyapunov equations as in \eqref{eq:LYAP1} and \eqref{eq:LYAP2}.

\subsection[Gradients of the modified objective]{Gradients of \eqref{eq:opt_modified}}
Let us derive the gradients of $f$ versus $\hat{A}$, $\hat{B}$, and $\hat{C}$. We define a gradient as follows.
\begin{definition}
[\!\!\!{\cite{van2008h2}}]
\label{def:gradient_model_base}
    The gradient of a real scalar smooth function $f(X)$ of a real matrix variable $X\in \Real^{n\times p}$ is the real matrix $\nabla_{X}f(X)\in \Real^{n\times p}$ defined by 
    \begin{align*}
        [\nabla_{X} f(X)]_{i,j} = \frac{\partial}{\partial X_{i,j}} f(X), \:\: i = 1,\ldots,n, \:\:j = 1,\ldots,p.
    \end{align*}
\end{definition}
\begin{prop}[\!{\cite{van2008h2}\cite{bunse2010h2}}]
\label{prop:garadient_model_base}
    The gradients $\nabla_{\hat{A}}f$, $\nabla_{\hat{B}}f$, and $\nabla_{\hat{C}}f$ of $f$ are given by
    \begin{align}
        &\nabla_{\hat{A}}f(\hat{A},\hat{B},\hat{C}) = 2(Q\hat{A}P + S^{\top}AR), \label{eq:gradient_Ahat}\\
        &\nabla_{\hat{B}}f(\hat{A},\hat{B},\hat{C}) = 2(S^{\top}B + Q\hat{B}), \label{eq:gradient_Bhat}\\
        &\nabla_{\hat{C}}f(\hat{A},\hat{B},\hat{C}) = 2(\hat{C}P-R),\label{eq:gradient_Chat}
    \end{align}
    where $P$, $Q$, $R$, and $S$ are the solutions to \eqref{eq:lyap1}, \eqref{eq:lyap2}, \eqref{eq:sylve1}, and \eqref{eq:sylve2}, respectively.
\end{prop}
\tblack{Since the set of stable matrices is open~\cite{orbandexivry2013nearest}, a stationary point can be defined as follows.}
\begin{definition}
    \tblack{A feasible triple $(\hat{A},\hat{B},\hat{C})$ for \eqref{eq:opt_modified} is called a stationary point if the gradients $(\nabla_{\hat{A}}f, \nabla_{\hat{B}}f, \nabla_{\hat{C}}f)$ for Proposition~\ref{prop:garadient_model_base} are all zero.}
\end{definition}
When the gradients $(\nabla_{\hat{A}}f, \nabla_{\hat{B}}f, \nabla_{\hat{C}}f)$ are all zero at certain points, the stationary conditions derived in \cite{bunse2010h2}, i.e., Wilson's first order necessary conditions, are satisfied.
However, note that these points may not necessarily belong to the feasible set ${\mathcal E}_{\mathrm{s}}$.
In this study, we derive these gradients directly from the measurement data without using the constant matrices $A$ and $B$ of the system \eqref{eq:lti}, and propose a gradient-based optimization algorithm for~\eqref{eq:opt_modified}.


\section{Problem formulation} \label{sec:formulation}
We formulate the problem addressed in this study under the assumption that the matrices $A\in \Real^{n\times n}$ and $B\in \Real^{n\times m}$ for the original system~\eqref{eq:lti} are unknown.
Note that this is a different setting from Section~\ref{sec:pre}.
Instead, assume that $N$ measurement data sets obtained from the true system \eqref{eq:lti} are given by
\begin{align}
\label{eq:DATA}
&\begin{pmatrix}
  x_{1,1} & x_{1,2} & \dots  & x_{1,L} \\
  x_{2,1} & x_{2,2} & \dots  & x_{2,L} \\
  \vdots & \vdots & \ddots & \vdots \\
  x_{N,1} & x_{N,2} & \dots  & x_{N,L}
\end{pmatrix},
\:\:
\begin{pmatrix}
  u_{1,1} & u_{1,2} & \dots  & u_{1,L} \\
  u_{2,1} & u_{2,2} & \dots  & u_{2,L} \\
  \vdots & \vdots & \ddots & \vdots \\
  u_{N,1} & u_{N,2} & \dots  & u_{N,L}
\end{pmatrix}
\end{align}
where $L$ is a total number of measurement instances for the data sets $i=1,2,\ldots,N$.
For \eqref{eq:DATA}, $x_{i,k+1}$, $x_{i,k}$, and $u_{i,k}$ satisfy
\begin{align}\label{eq:state_equation}
    x_{i,k+1}=Ax_{i,k}+Bu_{i,k},
\end{align}
where $i=1,\ldots,N$ and $k=1,\ldots,L-1$.

\begin{remark}\label{rem:state-data-assumption}
    While the assumption that state data is directly observable is restrictive, this is identical to the standard DMD framework \cite{kutz2016dynamic}, which has been applied in various fields such as neuroscience \cite{fieseler_unsupervised_2020}, fluid dynamics \cite{schmid_dynamic_2010}, epidemiology \cite{mustavee_linear_2022}, and financial engineering \cite{kuttichira_stock_2017}.
\end{remark}

In this paper, we address the following problem.
\begin{problem}\label{prob:data_driven_h2}
    Develop an algorithm for solving the $h^2$ MOR problem described in \eqref{eq:opt_modified} in a data-driven manner. 
    Specifically, construct a stable ROM using the measurement data \eqref{eq:DATA}, under the condition that the matrices $A$ and $B$ in \eqref{eq:lti} are unknown.
\end{problem}
Note that the gradients as expressed in Proposition \ref{prop:garadient_model_base} are not directly applicable due to the unknown nature of matrices $A$ and $B$.

Problem \ref{prob:data_driven_h2} is important for the following reason:
we examine a system with transfer function $H_{\textrm{SI}}$, which is obtained by system identification.
By the triangle inequality,
\begin{align}\label{eq:error_between_SI_and_True}
    \|H-\hat{H}\|_{h^2}
    &\leq \|H-H_{\textrm{SI}}\|_{h^2} + \|H_{\textrm{SI}}-\hat{H}\|_{h^2}.
\end{align}
\tblack{When using classical system identification~\cite{Ljung1999, Katayama2005}, the error~$\|H_{\textrm{SI}}-\hat{H}\|_{h^2}$ is taken into account.}
\tblack{
However, if errors $\|H-H_{\textrm{SI}}\|_{h^2}$ arise from system identification, minimizing $\|H_{\textrm{SI}}-\hat{H}\|_{h^2}$ may still not result in satisfactory ROMs in the $h^2$-norm. 
Furthermore, if the system identified is unstable, $\|H_{\textrm{SI}}-\hat{H}\|_{h^2}$ cannot be defined and the $h^2$ MOR may not be performed.}
To overcome these problems, this paper proposes a MOR method to optimize the left-hand side of \eqref{eq:error_between_SI_and_True} directly from the measurement data.


\section{Data-driven \texorpdfstring{$h^2$}{h2} Model Reduction for Discrete-time Systems} \label{sec:main}
We describe a data-driven $h^{2}$ MOR for \eqref{eq:opt_modified}. 
By making some assumptions on the measurement data of the original system \eqref{eq:lti}, we show that ROMs can be obtained directly from the data when the matrices $A$ and $B$ of \eqref{eq:lti} are unknown.

\subsection{Data-driven Reconstruction of the State Equation}\label{subsec:data-driven-state-equation}
To solve the discrete-time Sylvester equation \eqref{eq:sylve1} in a data-driven manner, we represent the system's state equation \eqref{eq:state_equation} using the measurement data \eqref{eq:DATA}.
First, we collect the data in the matrices
\begin{equation}
\label{eq:data_matrix}
\begin{aligned}
    X := &\big( x_{1,1}, x_{1,2}, \ldots, x_{1,L-1}, x_{2,1}, x_{2,2}, \ldots, x_{2,L-1}, \\
    & \ldots, x_{N,1}, x_{N,2}, \ldots, x_{N,L-1} \big)^{\top} 
    \in \mathbb{R}^{(NL-N) \times n}, \\
    X_{\mathrm{p}} := &\big( x_{1,2}, x_{1,3}, \ldots, x_{1,L}, x_{2,2}, x_{2,3}, \ldots, x_{2,L}, \\
    & \ldots, x_{N,2}, x_{N,3}, \ldots, x_{N,L} \big)^{\top} 
    \in \mathbb{R}^{(NL-N) \times n}, \\
    U := &\big( u_{1,1}, u_{1,2}, \ldots, u_{1,L-1}, u_{2,1}, u_{2,2}, \ldots, u_{2,L-1}, \\
    & \ldots, u_{N,1}, u_{N,2}, \ldots, u_{N,L-1} \big)^{\top} 
    \in \mathbb{R}^{(NL-N) \times m},
\end{aligned}
\end{equation}
where the $i$th row of $X_{\mathrm{p}}$ corresponds to the $i$th row of $X$, shifted by one step according to the system's state equation.

Using \eqref{eq:state_equation} and the data matrices \eqref{eq:data_matrix},
\begin{align}\label{eq:state_equation_data}
X_{\mathrm{p}}=XA^{\top}+UB^{\top}.
\end{align}
Multiplying this equation by $X^{\top}$ from the right yields
\begin{align}\label{eq:state_equation_data_mod}
    X_{\mathrm{p}}X^{\top}=XA^{\top}X^{\top}+UB^{\top}X^{\top}.
\end{align}
Here \eqref{eq:state_equation_data_mod} can be rewritten as
\begin{align}
    X_{\mathrm{p}}X^{\top} &= XZ_{A}^{\top}+UZ_{B}^{\top}, \label{eq:relation_x_z_matrix}
\end{align}
where
\begin{align}
    Z_A := XA,\quad Z_B := XB. \label{eq:ZAZB}
\end{align}
Furthermore, transposing \eqref{eq:state_equation_data_mod} yields 
\begin{align}
    XU_{B}^{\top}&=XX_{\mathrm{p}}^{\top}-Z_{A}X^{\top}, \label{eq:Bunknown_UB}
\end{align}
where
\begin{align}
    U_B &:= UB^{\top}. \label{eq:UB}
\end{align}
Note that $Z_A$ corresponds to the matrix of state data for the dual system in \eqref{eq:state_equation}.
In general, the solutions $Z_A$, $Z_B$, and $U_B$ of \eqref{eq:relation_x_z_matrix} and \eqref{eq:Bunknown_UB} are not uniquely determined as in \eqref{eq:ZAZB} and \eqref{eq:UB}. 
We denote the solutions obtained by solving \eqref{eq:relation_x_z_matrix} and \eqref{eq:Bunknown_UB} as $Z_{A,\text{data}}$,  $Z_{B,\text{data}}$, and $U_{B,\text{data}}$, respectively. 
To prove uniqueness, we assume the following.
\begin{assumption}\label{Assumptions}
Conditions for the proposed method:
\begin{enumerate}[label=(A\arabic*), wide, labelindent=0pt]
    \item The spectra of $X^{\dagger}Z_{A,\text{data}}$ and $\hat{A}^{-1}$ are disjoint. \label{assump_A1}
    \item The spectra of $X^{\dagger}(X_{\mathrm{p}}-U_{B,\text{data}})$ and $\hat{A}^{-1}$ are disjoint. \label{assump_A2}
    \item The eigenvalues of $\hat{A}$ are inside the unit circle and strictly positive, i.e., $\forall i\in \{1,2,\ldots,r\}, 0< |\lambda_{i}(\hat{A})| < 1$. \label{assump_A3} \vspace{-10pt}
    \item $\rank\!\begin{pmatrix} X & U \end{pmatrix}=n+m$. \label{assump_A4}
    \item $\rank X=n$. \label{assump_A5}
    \item $\rank U=m$. \label{assump_A6}
\end{enumerate}
\end{assumption}
\ref{assump_A1}-\ref{assump_A6} may be easily achieved for the following reasons.
\ref{assump_A1}-\ref{assump_A3} are conditions related to the eigenvalues of $\hat{A}$. As described in Section \ref{sec:algo}, the proposed algorithm allows $\hat{A}$ to be constructed to satisfy these assumptions. For details, refer to Section \ref{sec:algo}.

\ref{assump_A4}-\ref{assump_A6} are specifically related to the measurement data.
If the number of data sets \( N \) and measurement instances \( L \) satisfy $NL-N \geq n+m$, and the vectors \( (x_{i,k}, u_{i,k})^{\top} \), \( x_{i,k} \), and \( u_{i,k} \) are linearly independent for each \( i = 1, \ldots, N \) and \( k = 1, \ldots, L \), respectively, then the assumptions hold.
Furthermore, these assumptions are related to the persistency of excitation condition of the input \cite{willems_note_2005}. Specifically, when the input $u$ possesses sufficient persistent excitation, the Hankel matrix of $U$ becomes full rank, and the rank of $X$ coincides with the state dimension $n$. 
Therefore, for example, as shown in \cite{markovsky2023persistency}, if an input \( u \) that satisfies the persistency of excitation condition is generated, both $X$ and $U$ are guaranteed to be full rank.

\subsection{Construction of Gradient-like Matrices from Measurement Data}
In this subsection, we show that the gradient-like matrices corresponding to Proposition \ref{prop:garadient_model_base} can be derived from the data.
To obtain the gradient-like matrices, we solve the discrete-time Sylvester equations \eqref{eq:sylve1} and \eqref{eq:sylve2} in a data-driven manner, following \cite{banno2021data}\cite{vrabie2009adaptive}.
Multiplying \eqref{eq:sylve1} and \eqref{eq:sylve2} by the data matrix $X$ from the left, from \eqref{eq:state_equation_data}, \eqref{eq:ZAZB}, and \eqref{eq:UB}, we obtain
\begin{align}
    &XAR\hat{A}^{\top} + XB\hat{B}^{\top} = XR \notag\\
    &\Leftrightarrow Z_{A}R\hat{A}^{\top} + Z_{B}\hat{B}^{\top} = XR, \label{eq:sylve1_data_Bunknown}\\
    &XA^{\top}S\hat{A} - X\hat{C} = XS \notag\\
    &\Leftrightarrow (X_{\mathrm{p}}-U_{B})S\hat{A}-X\hat{C} = XS. \label{eq:sylve2_data_Bunknown}
\end{align}
Furthermore, multiplying $X^{\dagger}$ from the left-hand side of \eqref{eq:sylve1_data_Bunknown} and \eqref{eq:sylve2_data_Bunknown}, they can be rewritten as
\begin{align}
    &X^{\dagger}Z_{A}R\hat{A}^{\top} + X^{\dagger}Z_{B}\hat{B}^{\top} = X^{\dagger}XR, \label{eq:appendix_discrete_sylve1_lem2} \\
    &X^{\dagger}(X_{\mathrm{p}}-U_{B})S\hat{A}-X^{\dagger}X\hat{C} = X^{\dagger}XS, \label{eq:appendix_discrete_sylve2_lem2}
\end{align}
where the solutions $R$ and $S$ of \eqref{eq:appendix_discrete_sylve1_lem2} and \eqref{eq:appendix_discrete_sylve2_lem2} are denoted by $R_{\text{data}}$ and $S_{\text{data}}$, respectively.

\begin{lemma}
\label{lem:uniqueness_for_data_driven_equation_Bunknown}
    Let $Z_{A,\text{data}}=Z_{A}$, $Z_{B,\text{data}}=Z_{B}$, and $U_{B,\text{data}}=U_{B}$, where $Z_{A,\text{data}}$, $Z_{B,\text{data}}$, and $U_{B,\text{data}}$ are the solutions of \eqref{eq:ZAZB} and \eqref{eq:UB}, respectively.
    Suppose that \ref{assump_A1}-\ref{assump_A3}, and \ref{assump_A5} hold.
    Then, the solutions $R$ and $S$ to \eqref{eq:appendix_discrete_sylve1_lem2} and \eqref{eq:appendix_discrete_sylve2_lem2} exist uniquely.
\end{lemma}
\begin{proof}
    It is shown using the unique solution existence theorem for the Sylvester equation~\cite{simoncini2016computational}. See Appendix.
\end{proof}

\tblack{Using the relation $AR = (R - B \hat{B}^{\top}) (\hat{A}^{\dagger})^{\top}$ from~\eqref{eq:sylve1}, and the solutions $R_{\text{data}}$ of~\eqref{eq:appendix_discrete_sylve1_lem2} and $S_{\text{data}}$ of~\eqref{eq:appendix_discrete_sylve2_lem2},} the gradient-like matrices for~\eqref{eq:opt_modified} are calculated as 
\begin{align}
    &\tilde{\nabla}_{\hat{A}}f(\hat{A},\hat{B},\hat{C}) := 2(Q\hat{A}P + (S_{\text{data}}^{\top}R_{\text{data}}-S_{B,\text{data}}^{\top}\hat{B}^{\top})(\hat{A}^{\dagger})^{\top}), \label{eq:gradient_Ahat_data_Bunknown}\\
    &\tilde{\nabla}_{\hat{B}}f(\hat{A},\hat{B},\hat{C}) := 2(S_{B,\text{data}}^{\top} + Q\hat{B}), \label{eq:gradient_Bhat_data_Bunknown}\\
    &\tilde{\nabla}_{\hat{C}}f(\hat{A},\hat{B},\hat{C}) := 2(\hat{C}P-R_{\text{data}}),\label{eq:gradient_Chat_data_Bunknown}
\end{align}
\tblack{
where $S_{B,\text{data}}$ denotes the solution $S_{B}$ to the equation 
\begin{align}
\label{eq:Bunknown_BY}
    US_{B} = U_{B}S_{\text{data}},
\end{align}
and, as stated in Lemma~\ref{lem:SB}, it coincides with $B^{\top} S_{\text{data}}$ under certain conditions.}
In general, however, \eqref{eq:gradient_Ahat_data_Bunknown}-\eqref{eq:gradient_Chat_data_Bunknown} do not coincide with \eqref{eq:gradient_Ahat}-\eqref{eq:gradient_Chat} because $R_{\text{data}}$ and $S_{\text{data}}$ are different from the solutions $R$ and $S$ in~\eqref{eq:sylve1} and \eqref{eq:sylve2}.
Using such inaccurate gradients could lead to convergence to incorrect solutions. 
In what follows, we derive sufficient conditions on the data for the gradient-like matrices to coincide with the gradients of Proposition \ref{prop:garadient_model_base}.

\subsection{Data-driven \texorpdfstring{$h^2$}{h2} Model Reduction}
\label{subsec:model_reduction_Bunknown}
This subsection describes the data-driven model reduction for the case where $A\in \Real^{n\times n}$ and $B\in \Real^{n\times m}$ are unknown. 
To solve \eqref{eq:opt_modified} numerically, we compute the gradients of Proposition~\ref{prop:garadient_model_base} in a data-driven manner.
We first introduce two lemmas, whose proofs are in the appendix.

\begin{lemma}\label{lem:data_matrix_z}
    Let $NL-N\geq n+m$.
    Suppose that \ref{assump_A4}-\ref{assump_A6} hold. 
    Then, \eqref{eq:relation_x_z_matrix} has the unique solution $(Z_{A,\text{data}}, Z_{B,\text{data}}) = (Z_{A}, Z_{B})$.
\end{lemma}
Lemma \ref{lem:data_matrix_z} means that under assumptions, $Z_{A}$ and $Z_{B}$ given by~\eqref{eq:ZAZB} can be expressed using input and state data of system \eqref{eq:lti}.
In other words, $Z_{A}$ and $Z_{B}$ can be expressed as $Z_{A,\text{data}}$ and $Z_{B,\text{data}}$ without using the unknown matrices $A$ and $B$.
\begin{lemma}\label{lem:SB}
    Let $NL-N\geq n+m$.
    Suppose that \ref{assump_A4}-\ref{assump_A6} hold.
    Then, \eqref{eq:Bunknown_UB} and \eqref{eq:Bunknown_BY} have the unique solutions \eqref{eq:UB} and
    \begin{align}\label{eq:SB}
        S_{B} = B^{\top}S_{\text{data}},
    \end{align}
    respectively.
\end{lemma}
The following theorem presents the conditions under which the gradient-like matrices, obtained through a data-driven approach, become the gradients for \eqref{eq:opt_modified}.
\begin{theorem}
\label{thm:gradient_Bunknown}
    Let $NL-N\geq n+m$.
    Suppose that \ref{assump_A1}-\ref{assump_A6} hold. 
    Then, the gradient-like matrices of \eqref{eq:gradient_Ahat_data_Bunknown}-\eqref{eq:gradient_Chat_data_Bunknown} coincide with the gradients for \eqref{eq:opt_modified} specified in Proposition \ref{prop:garadient_model_base}, that is, $\tilde{\nabla}_{\hat{A}}f =\nabla_{\hat{A}}f$, $\tilde{\nabla}_{\hat{B}}f =\nabla_{\hat{B}}f$, and $\tilde{\nabla}_{\hat{C}}f =\nabla_{\hat{C}}f$ hold.
\end{theorem}
\begin{proof}
    It follows from Lemmas \ref{lem:uniqueness_for_data_driven_equation_Bunknown}, \ref{lem:data_matrix_z}, and \ref{lem:SB}.
    See Appendix.
\end{proof}
Theorem~\ref{thm:gradient_Bunknown} states that the gradients for \eqref{eq:opt_modified} are computed based on the data given by \eqref{eq:data_matrix} under \ref{assump_A1}-\ref{assump_A6}.
In fact, $R_{\text{data}}$, $S_{\text{data}}$, and $S_{B,\text{data}}$ in \eqref{eq:gradient_Ahat_data_Bunknown}, \eqref{eq:gradient_Bhat_data_Bunknown}, and \eqref{eq:gradient_Chat_data_Bunknown} can be calculated by solving \eqref{eq:appendix_discrete_sylve1_lem2}, \eqref{eq:appendix_discrete_sylve2_lem2}, and \eqref{eq:Bunknown_BY}.
Moreover, $(Z_A,Z_{B})$ and $U_{B}$ in \eqref{eq:appendix_discrete_sylve1_lem2}, \eqref{eq:appendix_discrete_sylve2_lem2}, and \eqref{eq:Bunknown_BY} can be calculated by solving \eqref{eq:relation_x_z_matrix} and \eqref{eq:Bunknown_UB}.
As mentioned earlier, if \ref{assump_A1}-\ref{assump_A6} are not satisfied, it is possible that the gradients obtained from the data differ from the true gradients as shown in Proposition \ref{prop:garadient_model_base}. 
When \ref{assump_A1}-\ref{assump_A6} regarding the data are satisfied and the gradients from Theorem~\ref{thm:gradient_Bunknown} are all zero, the stationary point condition described in Section~\ref{sec:pre} is satisfied. This fact becomes crucial in the construction of the proposed algorithm, as discussed in Section \ref{sec:algo}.

\tblack{
\begin{remark}\label{rem:C_identity}
Our framework assumes that the full state is directly measured, that is, \(p = n\) and, without loss of generality, \(C = I_n\).
Under this setting we have \(y_{i,k} = x_{i,k}\), and therefore the state sequence is exactly known and Theorem~\ref{thm:gradient_Bunknown} provides the true gradients of \eqref{eq:opt_modified}.
If only a part of the state is observed (\(p < n\)) or, more generally, if \(C \neq I_n\), the state \(x_{i,k}\) must be reconstructed, which inevitably introduces approximation errors.
In such a partially observed setting the matrices \(\tilde{\nabla}_{\hat{A}}f\), \(\tilde{\nabla}_{\hat{B}}f\), and
\(\tilde{\nabla}_{\hat{C}}f\) no longer coincide with the true gradients. Future work will address this case.
\end{remark}
}

\subsection[Gradients of the modified objective (A unknown, B known)]
{Gradients of \eqref{eq:opt_modified} for the case that $A$ is unknown and $B$ is known}
We derive the gradients for the case where the matrix $B$ for \eqref{eq:lti} is known.
Unlike the setting in Subsection \ref{subsec:model_reduction_Bunknown}, the conditions for the gradient-like matrices to coincide with the gradients of Proposition \ref{prop:garadient_model_base} are different.

\begin{corollary}
    \label{thm:gradient_Bknown}
    Let $NL-N\geq n$.
    Suppose that \ref{assump_A1}-\ref{assump_A3}, and \ref{assump_A5} hold. 
    Then, the gradient-like matrices of \eqref{eq:gradient_Ahat_data_Bunknown}-\eqref{eq:gradient_Chat_data_Bunknown} coincide with the gradients for \eqref{eq:opt_modified} specified in Proposition \ref{prop:garadient_model_base}, that is, $\tilde{\nabla}_{\hat{A}}f =\nabla_{\hat{A}}f$, $\tilde{\nabla}_{\hat{B}}f =\nabla_{\hat{B}}f$, and $\tilde{\nabla}_{\hat{C}}f =\nabla_{\hat{C}}f$ hold true.
\end{corollary}
\begin{proof}
    The proof follows directly from that of Theorem~\ref{thm:gradient_Bunknown}.
\end{proof}

The next corollary is the condition under which the gradients of Theorem~\ref{thm:gradient_Bunknown} and Corollary~\ref{thm:gradient_Bknown} coincide.
\begin{corollary}\label{cor:coincide_condition}
    Under \ref{assump_A1}-\ref{assump_A6}, the gradients \eqref{eq:gradient_Ahat_data_Bunknown}-\eqref{eq:gradient_Chat_data_Bunknown} with unknown $B\in \Real^{n\times m}$ coincide with the gradients with known $B\in \Real^{n\times m}$.
\end{corollary}
\begin{proof}
    From Lemma \ref{lem:SB}, we obtain
    $\nabla_{\hat{A}}f
    = 2(Q\hat{A}P + (S^{\top}R-S_{B}^{\top}\hat{B}^{\top})(\hat{A}^{-1})^{\top})
    = 2(Q\hat{A}P + S^{\top}(R-B\hat{B}^{\top})(\hat{A}^{-1})^{\top})$,
    $\nabla_{\hat{B}}f= 2(S_{B}^{\top} + Q\hat{B}) = 2(S^{\top}B + Q\hat{B})$.
\end{proof}

\section{Gradient-based data-driven \texorpdfstring{$h^2$}{h2} MOR algorithm} 
\label{sec:algo}
In this section, we present the Algorithm~\ref{alg:Bunknown} which solves the problem \eqref{eq:opt_modified} when $A\in \Real^{n\times n}$ and $B\in \Real^{n\times m}$ are unknown. 
\tblack{Consider the set $\mathcal E:=\{\hat{\theta}=(\hat A,\hat B,\hat C)\mid 0<|\lambda_i(\hat{A})|<1,\ i=1,\dots,r\}$ which satisfy~\ref{assump_A3}.}
\tblack{Since $\mathcal E$ is open (proved in Proposition~\ref{prop:existence_of_alpha} below), if the gradients from Theorem~\ref{thm:gradient_Bunknown} are all zero at a point in $\mathcal E$, that point is a stationary point.}

We briefly describe Algorithm~\ref{alg:Bunknown}.
First, we generate the initial reduced matrices $\hat{\theta}_{1}=(\hat{A}_{(1)}, \hat{B}_{(1)}, \hat{C}_{(1)})\in \mathcal E$, satisfying~\ref{assump_A3}.
In step 1 and 2 of Algorithm~\ref{alg:Bunknown}, $Z_{A,\text{data}}$, $Z_{B,\text{data}}$, and $U_{B,\text{data}}$, which are needed to solve the discrete-time Sylvester equations \eqref{eq:appendix_discrete_sylve1_lem2} and \eqref{eq:appendix_discrete_sylve2_lem2} with $(Z_{A},Z_B,U_B)=(Z_{A,\text{data}},Z_{B,\text{data}},U_{B,\text{data}})$, are generated. 
At each iteration~$\ell$, after computing the gradients, the reduced matrices $\hat{\theta}_{\ell}:=(\hat A_{(\ell)},\hat B_{(\ell)},\hat C_{(\ell)})$
are updated so that the objective function for \eqref{eq:opt_modified} becomes smaller.
In the while statement from step 10 to 16, the backtracking method is executed. 
\begin{prop}\label{prop:existence_of_alpha}
\tblack{
Suppose $\hat \theta_{\ell}\in\mathcal E$. 
Then there exists $\alpha^\star_{\ell}>0$ such that for all $\alpha\in(0,\alpha^\star_{\ell})$,}
\begin{align*}
    \tblack{\bar\theta:=\hat\theta_{\ell}-\alpha\,\tilde{\nabla}f(\hat \theta_{\ell})\in\mathcal E.}
\end{align*}
\end{prop}
\begin{proof}
\tblack{See Appendix.}
\end{proof}

The step-size $\alpha$ can be determined by the backtracking method such that the Armijo rule~$f(\hat\theta_{\ell + 1})\leq f(\hat{\theta}_{\ell}) - c_1 \alpha_{\ell}\:\|\tilde \nabla f(\hat{\theta}_{\ell})\|^{2}$
and \ref{assump_A1}-\ref{assump_A3} are satisfied.
\tblack{From a numerical perspective, by continuity of eigenvalues, the spectrum of $\hat{A}^{-1}$ does not intersect with that of $X^{\dagger}Z_{A,\text{data}}$ nor with that of $X^{\dagger}(X_{\mathrm{p}}-U_{B,\text{data}})$, thus~\ref{assump_A1} and \ref{assump_A2} are readily satisfied.
}
\tblack{Furthermore, from 
Theorem~\ref{thm:gradient_Bunknown}, 
Proposition~\ref{prop:existence_of_alpha}, and the smoothness of $f$, there exists $\alpha>0$ satisfying the Armijo rule and preserving feasibility.}

The convergence of Algorithm~\ref{alg:Bunknown} is guaranteed by the following theorem. Note that \ref{assump_A1}-\ref{assump_A3} are satisfied in each iteration.
\begin{theorem}\label{thm:convergence}
    \tblack{Assume \ref{assump_A4}–\ref{assump_A6} and that the sequence $\{\hat\theta_\ell\}$ with $\hat\theta_\ell=(\hat A_{(\ell)},\hat B_{(\ell)},\hat C_{(\ell)})$ generated by Algorithm~\ref{alg:Bunknown} with $tol=0$ is bounded and $\lim_{\ell\rightarrow\infty}\hat A_{(\ell)}$ satisfies \ref{assump_A3}. Then $\{\hat\theta_\ell\}$ converges to a stationary point of~\eqref{eq:opt_modified}.}
\end{theorem}
\begin{proof}
\tblack{
It follows from~\cite[Thm. 3.2]{attouch2013convergence}.
See Appendix.
}
\end{proof}

\tblack{Note that while the existence of a stationary point for \eqref{eq:opt_modified} is non-trivial, it is guaranteed when the assumptions of Theorem~\ref{thm:convergence} are satisfied.}
Next, we provide the computational complexity analysis of Algorithm~\ref{alg:Bunknown}.
\begin{theorem}\label{thm:cost}
Assume that in Algorithm~\ref{alg:Bunknown} the backtracking terminates in at most $J$ trials per iteration, and the algorithm performs $\ell_{\text{max}}$ iterations in total.
Then, the overall computational complexity of Algorithm~\ref{alg:Bunknown} is given by 
\begin{align*}
    \mathcal{O}((NL-N)^2 (n+m) + ((NL-N)n^2+ n^3+r^3) \:J\:\ell_{\text{max}}).
\end{align*}
\end{theorem}
\begin{proof}
    See Appendix.
\end{proof}
\tblack{
A widely–used alternative for constructing a ROM directly from snapshot data is DMDc \cite{proctor2016dynamic, kutz2016dynamic}.
Given the snapshot matrix $X$, $X_{\text{p}}$, 
$U$ by \eqref{eq:data_matrix}, the computational complexity by DMDc~\cite{proctor2016dynamic, kutz2016dynamic} is $\mathcal{O}\bigl((NL-N)^{2}(r+m)+(NL-N)nr\bigr)$.
Although Theorem~\ref{thm:cost} indicates that the proposed method incurs a higher computational complexity than DMDc, it can generate a stable ROM with a smaller $h^{2}$ error, as described in Section~\ref{subsec:computational_neuroscience}.
}

\begin{remark}
   Since \eqref{eq:opt_modified} is a nonconvex optimization problem, the choice of the initial point for Algorithm~\ref{alg:Bunknown} is important.  
   As the initial point, it is possible to use ROMs obtained via existing data-driven MOR methods, such as DMDc~\cite{kutz2016dynamic}\cite{proctor2016dynamic}, 
   Loewner framework~\cite{mayo2007framework,gosea2021data}, and 
   data-driven BT~\cite{gosea2022data,burohman2023data,ma2011reduced}.
   Such initialization often produces high-quality ROMs in the $h^2$-norm.
   However, these methods require different types of data, and therefore, they may not be directly applicable in some cases.
\end{remark}

\begin{remark}\label{rem:model_base_data_base_alg}
\tblack{
Consider system~\eqref{eq:state_equation} with known \(A,B\), and compare Algorithm~\ref{alg:Bunknown} (data-driven) with its model-based counterpart, the model-based $h^2$ MOR using gradients.
Fix \(\alpha_{\mathrm{init}}, c_1, \rho, \mathit{tol}, \hat{\theta}^{(1)}\).
Suppose that the data $(x_{i,k}$, $u_{i,k})$, $(i=1,2,\ldots,N, k=1,2,\ldots,L)$ are generated from~\eqref{eq:state_equation} and \ref{assump_A1}-\ref{assump_A6} hold.
Then, by Theorem~\ref{thm:gradient_Bunknown}, \(\nabla f\) and \(\tilde{\nabla} f\) coincide for the same $\hat{\theta}=(\hat{A},\hat{B}, \hat{C})$, and the \(\alpha_\ell\) determined by the backtracking method also coincide.
Therefore, the sequences $\{\hat\theta_{\,\ell}\}$ generated by the two methods are identical.}
\end{remark}

\begin{figure}[!t]
\begin{algorithm}[H]
    \caption{Solve \eqref{eq:opt_modified} for the case that $A$ and $B$ are unknown}
    \label{alg:Bunknown}
    \begin{algorithmic}[1]
    \REQUIRE Measurement data $(x_{i,k}$, $u_{i,k})$, $(i=1,2,\ldots,N, k=1,2,\ldots,L)$, dimension $r$ for \eqref{eq:reduced_lti}, initial step-size $\alpha_\text{init}$, Armijo parameter $c_1$, search control parameter $\rho$, tolerance $tol$, the initial real matrices $\hat{\theta}_{1}$ which satisfy~\ref{assump_A3}
    \ENSURE $\hat{\theta}=(\hat{A}, \hat{B}, \hat{C})\in \Real^{r\times r}\times \Real^{r\times m}\times \Real^{n\times r}$
    \STATE Generate $Z_{A,\text{data}}$ and $Z_{B,\text{data}}$ from \eqref{eq:relation_x_z_matrix} with $(Z_{A},Z_{B})=(Z_{A,\text{data}},Z_{B,\text{data}})$
    \STATE Generate $U_{B,\text{data}}$ from \eqref{eq:Bunknown_UB} with $(Z_{A},U_{B})=(Z_{A,\text{data}},U_{B,\text{data}})$
    \FOR{$\ell=1,2,\ldots$}
    \STATE Solve \eqref{eq:lyap1} and \eqref{eq:lyap2} for $P$ and $Q$ with $\hat{\theta}=\hat{\theta}_{\ell}$ 
    \STATE Solve \eqref{eq:appendix_discrete_sylve1_lem2} and \eqref{eq:appendix_discrete_sylve2_lem2} for $R$ and $S$ with $\hat{\theta}=\hat{\theta}_{\ell}$ 
    \STATE Calculate gradient-like matrices $\tilde{\nabla}_{\hat{A}}f(\hat{\theta}_{\ell})$, $\tilde{\nabla}_{\hat{B}}f(\hat{\theta}_{\ell})$, and $\tilde{\nabla}_{\hat{C}}f(\hat{\theta}_{\ell})$ with $S_{B,\text{data}}$ obtained by \eqref{eq:Bunknown_BY}
    \STATE Define $\tilde{\nabla} f(\hat{\theta}_{\ell}):=(\tilde{\nabla}_{\hat{A}}f(\hat{\theta}_{\ell}), \tilde{\nabla}_{\hat{B}}f(\hat{\theta}_{\ell}), \tilde{\nabla}_{\hat{C}}f(\hat{\theta}_{\ell}))$
    \STATE \textbf{if} $\|\tilde{\nabla} f(\hat{\theta}_{\ell})\|_{F}^{2} < tol$ \textbf{then break}
    \STATE $\alpha_{\ell} = \alpha_\text{init}$
    \WHILE{true}
    \STATE 
    $\bar{\theta}:=\hat{\theta}_{\ell}-\alpha_{\ell}\tilde{\nabla} f(\hat{\theta}_{\ell})$
    \STATE Solve \eqref{eq:lyap1} for $P$ with $\hat{\theta}=\bar{\theta}$
    \STATE Solve \eqref{eq:appendix_discrete_sylve1_lem2} for $R$ with $\hat{\theta}=\bar{\theta}$
    \STATE \textbf{if} $f(\bar{\theta})\leq f(\hat{\theta}_{\ell}) - c_{1}\alpha_{\ell}\:\|\tilde{\nabla} f(\hat{\theta}_{\ell})\|_{F}^{2}$ and $\forall i\:, \:0<|\lambda_{i}(\bar{A})|< 1$ \textbf{then} $\hat{\theta}_{\ell+1}:=\bar{\theta}$ \textbf{break}
    \STATE $\alpha_{\ell}\leftarrow\rho \alpha_{\ell}$
    \ENDWHILE
    \ENDFOR
    \end{algorithmic}
\end{algorithm}
\end{figure}

\section{Application to Computational Neuroscience} \label{subsec:computational_neuroscience}
In Computational Neuroscience, it is crucial to elucidate how the neuronal network of Caenorhabditis elegans (C. elegans) regulates behavioral states \cite{fieseler_unsupervised_2020,kato_global_2015, nichols_global_2017}.
Specifically, the behavioral states of C. elegans can be analyzed by observing the fluorescence intensity of neuronal groups \cite{kato_global_2015, nichols_global_2017}.
Moreover, it is known that changes in oxygen concentration affect the neural activity of C. elegans  \cite{kato_global_2015, nichols_global_2017}.

As in \cite{fieseler_unsupervised_2020}, we model the relationship between the neural activity of C. elegans, regarded as the state, and changes in oxygen concentration, regarded as the control input, using the following system with noise $\epsilon_{k}$:
\begin{empheq}[left= \empheqlbrace]{align} \label{eq:lti_noise}
    \begin{aligned}
    \tilde{x}_{k+1} &=A\tilde{x}_{k} + Bu_{k}, \\
    x_{k} &= \tilde{x}_{k} + \beta \epsilon_{k},
    \end{aligned}
\end{empheq}
where $k=1,\ldots,L-1$, and the noise \( \epsilon_{k} \) was generated using the MATLAB~\lstinline[style=Matlab-editor]{randn}~command.
Furthermore, we demonstrate that the ROMs obtained by the proposed data-driven method can reconstruct the original system more accurately in the sense of the $h^2$-norm compared to DMDc \cite{proctor2016dynamic, kutz2016dynamic}.

In experiments, \tblack{we fixed the seed with \lstinline[style=Matlab-editor]{rng(0,'twister')} to ensure exact reproducibility.}
Furthermore, the initial step-size $\alpha_\text{init}$, Armijo parameter $c_1$, search control parameter $\rho$ required for the backtracking method, and tolerance $tol$ in Algorithm~\ref{alg:Bunknown} were $\alpha=1$, $c=10^{-4}$, $\rho=0.5$, and $tol=10^{-3}$ respectively.

\subsection{Problem setting} \label{subsubsec:problem_set}
First, we describe the measurement data for the original system \eqref{eq:lti_noise}.
The data utilized in this study are accessible at \url{https://osf.io/a64uz/} and are linked to two prior experimental papers \cite{nichols_global_2017, fieseler_unsupervised_2020}.
Using the data corresponding to \cite{nichols_global_2017}, we constructed the state data \( X \) and \( X_\mathrm{p} \) with noise and \( \tilde{X} \) and \( \tilde{X}_\mathrm{p} \) without noise, as well as the control input data \( U \), as defined by \eqref{eq:data_matrix} and \eqref{eq:lti_noise}. 
More specifically, we utilized the calcium imaging data from the \textit{npr-1} mutant strain, which is sensitive to changes in oxygen concentration, during the pre-lethargus stage of this strain as the state data. 
The changes in oxygen concentration were used as the control input. 
Note that, in this dataset, the parameters for \eqref{eq:data_matrix} were set as \( N = 1 \), \( L = 4045 \), \( n = 114 \), and \( m = 1 \).

To incorporate broader temporal information, we consider a time delay of \( h_x = 9 \) steps for the state data \( X \) and \( h_u = 2 \) steps for the input data \( U \). Accordingly, we construct the following Hankel matrices \(X^{(h_x)}\) from \( X \) in \eqref{eq:data_matrix} with $\bar{L}:=L-h_{x}$:
\begin{align*}
X^{(h_x)} = 
\begin{bmatrix}
x_1 & x_2 & \cdots & x_{\bar{L}} \\
x_2 & x_3 & \cdots & x_{\bar{L} + 1} \\
\vdots & \vdots & \ddots & \vdots \\
x_{h_x} & x_{h_x+1} & \cdots & x_{\bar{L}+h_{x}-1}
\end{bmatrix}
\in \mathbb{R}^{h_x n\times \bar{L}}.
\end{align*}
In the similar way, $\tilde{X}^{(h_x)}\in \mathbb{R}^{h_x n\times \bar{L}}$, $X^{(h_x)}_{\mathrm{p}}\in \mathbb{R}^{h_x n\times \bar{L}}$, $\tilde{X}^{(h_x)}_{\mathrm{p}}\in \mathbb{R}^{h_x n\times \bar{L}}$, and $U^{(h_u)} \in \mathbb{R}^{h_u m\times \bar{L}}$ are constructed from $X$, $X_{\mathrm{p}}$, and $U$ in \eqref{eq:data_matrix}, respectively.
Therefore, we consider the following extended state-space model for \eqref{eq:lti_noise}:
\begin{empheq}[left= \empheqlbrace]{align} \label{eq:lti_noise_H}
    \begin{aligned}
    \tilde{x}_{k+1}^{(h_x)} &=A^{(h_x)}\tilde{x}_{k}^{(h_x)} + B^{(h_u)}u_{k}^{(h_u)}, \\
    x_{k}^{(h_x)} &= \tilde{x}_{k}^{(h_x)} + \beta \epsilon_{k}^{(h_x)}
    \end{aligned}
\end{empheq}
where \( x_{k}^{(h_x)} \), \( \tilde{x}_{k}^{(h_x)} \), and \( u_{k}^{(h_u)} \) represent the \( k \)th columns of the matrices \( X^{(h_x)} \), \( \tilde{X}^{(h_x)} \), and \( U^{(h_u)} \), respectively, and \( \epsilon_{k}^{(h_x)} = (\epsilon_{k}, \dots, \epsilon_{k+h_x -1}) \).
Using \eqref{eq:lti_noise_H}, $A^{(h_x)}$ and $B^{(h_u)}$ matrices are calculated from $X^{(h_x)}$, $X^{(h_x)}_{\mathrm{p}}$, and $U^{(h_u)}$ as follows:
\begin{align*}
    A^{(h_x)} &= 
    \left[\tilde{X}^{(h_x)}_{\mathrm{p}} 
    \begin{bmatrix}
        \tilde{X}^{(h_x)} \\ 
        U^{(h_u)}
    \end{bmatrix}^{\dagger}\right](:, 1:h_x n), \\
    B^{(h_u)} &= 
    \left[\tilde{X}^{(h_x)}_{\mathrm{p}} 
    \begin{bmatrix}
        \tilde{X}^{(h_x)} \\ 
        U^{(h_u)}
    \end{bmatrix}^{\dagger}\right](:, h_x n+1:h_x n + h_u m).
\end{align*}

To assess the effectiveness of the proposed method on the original system~\eqref{eq:lti_noise_H}, we report the relative error $\|H-\hat H\|_{h^2} / \|H\|_{h^2}$ where $H$ denotes the transfer functions of the original system~\eqref{eq:lti_noise_H}.
\tblack{
It should be noted that, data-driven model reduction methods (e.g., data-driven BT) that use frequency response data or impulse response data as input are difficult to apply directly in this context.
}

\subsection{Numerical results} \label{subsubsec:results}
We present the results of applying Algorithm~\ref{alg:Bunknown} to the given state and input data. The initial point for Algorithm~\ref{alg:Bunknown} is set to the reduced matrices obtained via DMDc, and \( X^{(h_x)}, U^{(h_u)}, X^{(h_x)}_{\mathrm{p}}, \hat{A} \) satisfy \ref{assump_A1}-\ref{assump_A6}.

\begin{figure}[htbp]
    \centering
    \includegraphics[width=8.25cm]{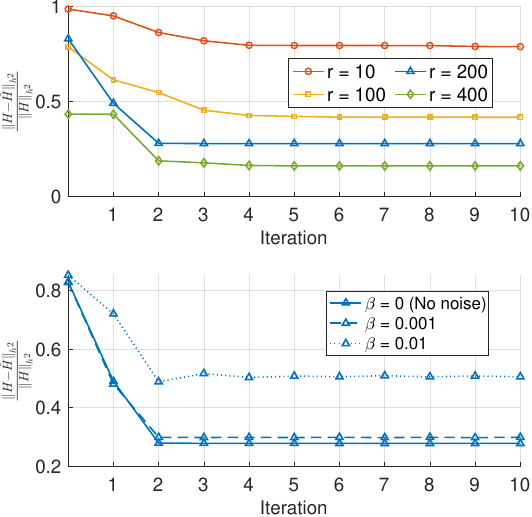}
    \caption{\tblack{Convergence behavior of \eqref{eq:opt_modified}: (upper) variation with $r$ ($\beta=0$, i.e., no noise), (lower) variation with $\beta$ ($r=200$).}}
    \label{fig:neuro_convergence}
\end{figure}

Fig.~\ref{fig:neuro_convergence} (upper) shows the \( h^2 \)-norm plotted for different values of \( r \) when the noise coefficient \( \beta \) is $0$. Similar to conventional model reduction methods, it is confirmed that larger values of \( r \) lead to ROMs that better approximate the original system~\eqref{eq:lti_noise_H}.
\tblack{Nevertheless, even at \(r = 400\) the relative \(h^{2}\)-error plateaus at about \(15\%\) because an intrinsically nonlinear model is being approximated by a purely linear ROM.}

Fig.~\ref{fig:neuro_convergence} (lower) illustrates the \( h^2 \)-norm for different noise coefficients \( \beta \) when \( r = 200 \). When data satisfying the assumptions is used, a better ROM in terms of the \( h^2 \)-norm is obtained compared to the initial point provided by DMDc. In the case of higher noise (\(\beta = 0.01\)), the noise impact cannot be ignored, causing the gradients to deviate significantly from the true gradients, which results in limited improvement. Conversely, when the noise level is small (\(\beta = 0.001\)), we observe a 50--60\% improvement, similar to the no-noise case (\(\beta = 0\)).

\section{Concluding remarks} \label{sec:conclusion}
In this study, we developed a data-driven model reduction method for linear discrete-time systems that constructs optimal ROMs in the $h^{2}$-norm sense. 
Through both theoretical and numerical analyses, we demonstrated that if the measurement data snapshots satisfy some assumptions, an $h^{2}$-optimal ROM can be obtained compared to the initial reduced model. 
Moreover, when modeling changes in neural activity of C. elegans, our approach was shown to outperform DMDc in terms of achieving an improved ROM under the $h^{2}$-norm sense. 

Future work includes: (i) more computationally efficient algorithms for large-scale problems; (ii) detailed evaluations in additional domains—preliminary results show our method also outperforms DMDc for daily infection modeling in epidemiology~\cite{mustavee_linear_2022}; and (iii) a Riemannian reformulation of \eqref{eq:opt_modified}, where \(\hat A\) is inherently stable and the explicit stability constraint becomes unnecessary~\cite{sato2019riemannian,obara2023stable}.

\section*{Acknowledgment}
This work was supported by JSPS KAKENHI under Grant Numbers 23K28369 and 25KJ0986.

\appendix
\label{appendix}
\begin{proof}[Proof of Lemma \ref{lem:uniqueness_for_data_driven_equation_Bunknown}] 
    Since~\ref{assump_A5}, $X^{\dagger}X=I_{n}$.
    From \cite{simoncini2016computational},
    under \ref{assump_A1} and \ref{assump_A3}, \eqref{eq:appendix_discrete_sylve1_lem2} has a unique solution $R_{\text{data}}=R$.
    Similarly, under \ref{assump_A2} and \ref{assump_A3}, \eqref{eq:appendix_discrete_sylve2_lem2} has a unique solution $S_{\text{data}}=S$.
\end{proof}

\begin{proof}[Proof of Lemma \ref{lem:data_matrix_z}]
    Since~\eqref{eq:ZAZB} are solutions of \eqref{eq:relation_x_z_matrix}, $\rank \!\begin{pmatrix} X & U \end{pmatrix}=\rank\!\begin{pmatrix} X & U & X_{\mathrm{p}}X^{\top} \end{pmatrix}$ holds.
    Thus, \ref{assump_A4} implies that \eqref{eq:relation_x_z_matrix} has the unique solution $(Z_{A,\text{data}}, Z_{B,\text{data}}) = (Z_{A}, Z_{B})$.
\end{proof}

\begin{proof}[Proof of Lemma \ref{lem:SB}]
    \eqref{eq:Bunknown_UB} has the unique solution \eqref{eq:UB} when \ref{assump_A5} holds.
    From \eqref{eq:Bunknown_BY}, we obtain $US_{B}=U_{B}S_{\text{data}}=UB^{\top}S_{\text{data}}$.
    Therefore, from~\ref{assump_A6}, $U^{\dagger}U=I_{m}$ holds. Thus, \eqref{eq:SB} holds. 
\end{proof}

\begin{proof}[Proof of Theorem~\ref{thm:gradient_Bunknown}]
    Since $Z_{A,\text{data}}=Z_{A}$,  $Z_{B,\text{data}}=Z_{B}$, and $U_{B,\text{data}}=U_{B}$ hold from Lemmas \ref{lem:data_matrix_z} and \ref{lem:SB}, the unique solutions \( R \) and \( S \) of \eqref{eq:sylve1} and \eqref{eq:sylve2}, respectively, coincide with the unique solutions \( R_{\text{data}} \) and \( S_{\text{data}} \) of \eqref{eq:appendix_discrete_sylve1_lem2} and \eqref{eq:appendix_discrete_sylve2_lem2}, as stated in Lemma~\ref{lem:uniqueness_for_data_driven_equation_Bunknown}.
    Furthermore, from Lemma \ref{lem:SB} and \eqref{eq:sylve1},
    we have $S^{\top}AR = S^{\top}(R-B\hat{B}^{\top})(\hat{A}^{-1})^{\top} = (S^{\top}R-S_{B}^{\top}\hat{B}^{\top})(\hat{A}^{-1})^{\top}$.
    Therefore, under \ref{assump_A1}-\ref{assump_A6}, the gradients of Proposition \ref{prop:garadient_model_base} are rewritten as \eqref{eq:gradient_Ahat_data_Bunknown}-\eqref{eq:gradient_Chat_data_Bunknown}.
\end{proof}

\begin{proof}[Proof of Proposition~\ref{prop:existence_of_alpha}]
    \tblack{
    The sets ${\mathcal E}_{\mathrm{s}}$ by \eqref{eq:feasible_set} (open; e.g.\ \cite{orbandexivry2013nearest}) and
    ${\mathcal E}_{\mathrm{n}}:=\{\hat\theta\mid\det \hat A\neq 0\}$ (open by continuity of $\det$) are open, hence 
    ${\mathcal E}=\mathcal E_{\mathrm{s}}\cap\mathcal E_{\mathrm{n}}$ is open. 
    Thus, there is $\delta>0$ with $\{ \bar\theta \mid \|\bar\theta - \hat{\theta}_{\ell}\|_F < \delta \}
    \subset\mathcal E$. 
    Therefore, when $\tilde\nabla f(\hat\theta_{\ell})\neq 0$, if $\alpha^\star_{\ell}:=\varepsilon/\|\tilde{\nabla}_{\hat A}f(\hat A_{(\ell)})\|_F$, then for any $\alpha\in (0, \alpha^\star_{\ell})$, $\bar\theta\in \mathcal E$ holds. Note that when $\tilde\nabla f(\hat\theta_{\ell})= 0$, $\bar\theta=\hat\theta$.}
\end{proof}

\begin{proof}[Proof of Theorem~\ref{thm:convergence}]
    \tblack{For the same $\hat{\theta}=(\hat A,\hat B,\hat C)$, the objective functions coincide; hence, by Theorem~\ref{thm:gradient_Bunknown}, $\nabla f=\tilde{\nabla} f$. Moreover, since both objectives are smooth, their Hessians also coincide, i.e., $\nabla^2 f=\tilde{\nabla}^2 f$.}

    \tblack{
    We define the closure
    $K:=\overline{\{\hat{\theta}_\ell\mid\ \ell\ge1\}}$ of the sequence generated by Algorithm~\ref{alg:Bunknown}.
    Since $\{\hat{\theta}_{\ell}\}$ is bounded and the limit $\lim_{\ell\to\infty}\hat A_{(\ell)}$ satisfies~\ref{assump_A3}, $K$ is a compact subset of ${\mathcal E}$. 
    Furthermore, since $f$ is smooth, $\nabla^2 f$ is continuous; hence, by the extreme value theorem,
    $L_{\nabla}:=\sup_{x\in K}\|\nabla^2 f(x)\|<\infty$.
    Therefore, 
    }
    \begin{align}\label{eq:lip}
    \tblack{
        \|\nabla f(x)-\nabla f(y)\|\le L_{\nabla}\|x-y\|
        \quad \forall\, x,y\in K.
    }
    \end{align}
    
    \tblack{
    Also, letting $s_{\ell}:=\hat{\theta}_{\ell+1}-\hat{\theta}_{\ell}=-\alpha_{\ell}\nabla f(\hat{\theta}_{\ell})$, we obtain}
    \begin{align}\label{eq:ABS_H1}
        \tblack{f(\hat\theta_{\ell + 1})\leq f(\hat{\theta}_{\ell}) - 
        \frac{c_1}{\alpha_{\ell}}\|s_{\ell}\|^{2}
        \leq f(\hat{\theta}_{\ell}) - 
        \frac{c_1}{\alpha_{\text{init}}}\|s_{\ell}\|^{2}.}
    \end{align}
    \tblack{
    From the Descent lemma~\cite[Lem.~3.1]{attouch2013convergence}, \eqref{eq:lip}, and the Armijo rule, we get
    $\alpha_{\ell}\leq 2(1-c_1)/ L_{\nabla}$.
    Furthermore, by the boundedness of $\{\hat{\theta}_\ell\}$ and the continuity of $\nabla f$, there exists $G<\infty$ such that $\|\nabla f(\hat{\theta}_\ell)\|\le G$. 
    Now, since $\lim_{\ell\rightarrow\infty}\hat A_{(\ell)}$ satisfies~\ref{assump_A3} and ${\mathcal E}$ is open,
    for any $x\in K$ and any $y\in \partial{\mathcal E} := \{\hat{\theta}\mid \min_i |\lambda_i(\hat A)|=0 \text{ or }\max_i |\lambda_i(\hat A)|=1\}$, there exists $\delta>0$ such that $\|x-y\|\ge \delta$. 
    Therefore, $\alpha_{\ell}\leq \delta/G\leq\delta/\|\nabla f(\hat{\theta}_{\ell})\|$ holds.
    Thus, for any $\ell\geq 1$, $\alpha_{\ell}\geq \underline{\alpha}:=\rho\cdot\min \{\alpha_{\text{init}}, 2(1-c_1)/L_{\nabla}, \delta/G\}$.
    As a result, we obtain}
    \begin{align}\label{eq:ABS_H2}
        \tblack{\|\nabla f(\hat{\theta}_{\ell + 1})\|\leq \|\nabla f(\hat{\theta}_{\ell})\| + L_{\nabla}\|s_\ell\|
        \leq (\underline{\alpha}^{-1}+ L_{\nabla})\|s_\ell\|,}
    \end{align}
    \tblack{where the first inequality follows from~\eqref{eq:lip}, and the second from the definition of $s_\ell$.}
    
    \tblack{
    Lastly, $\{\hat{\theta}_{\ell}\}$ is bounded and, by the Bolzano–Weierstrass theorem, there exists a convergent subsequence $\{\hat\theta_{\ell_{j}}\}$ whose limit lies in $K$. 
    Furthermore, since $f$ is real-analytic, it is a KL function.
    Hence, by \eqref{eq:ABS_H1}, \eqref{eq:ABS_H2}, and~\cite[Thm.~3.2]{attouch2013convergence}, Algorithm~\ref{alg:Bunknown} converges to a stationary point in $K$.
    }
\end{proof}

\begin{proof}[Proof of Theorem~\ref{thm:cost}]
    \tblack{To solve \eqref{eq:relation_x_z_matrix} for $(Z_{A},Z_{B})$, it is necessary to solve $NL-N$ linear equations of $(NL-N)\times (n+m)$. 
    Considering the least-squares problem, the complexity is $\mathcal{O}((NL-N)^2 (n+m))$. 
    Similarly, the complexity required to solve \eqref{eq:Bunknown_UB} for $U_{B}$ is $\mathcal{O}((NL-N)^2 n)$.
    The discrete‐time Lyapunov equations~\eqref{eq:lyap1} and \eqref{eq:lyap2} are each solvable in \(\mathcal O(r^{3})\) via the Bartels-Stewart method \cite{simoncini2016computational}.
    For the Sylvester equations~\eqref{eq:appendix_discrete_sylve1_lem2} and \eqref{eq:appendix_discrete_sylve2_lem2}, forming $X^{\dagger}$ and the matrix products
    $X^{\dagger}Z_{A}$, $X^{\dagger}Z_{B}\hat{B}^{\top}$, and $X^{\dagger}(X_{\mathrm{p}}-U_{B})$ each require
    $\mathcal{O}((NL-N)n^2)$.
    With \(X^{\dagger}X = I_n\), the resulting \(n\times n\) Sylvester equations are solved in \(\mathcal O(n^{3})\) via the Bartels-Stewart method \cite{simoncini2016computational}, giving a total cost $\mathcal{O}((NL-N)n^2+n^3)$.
    }
\end{proof}

\bibliographystyle{IEEEtran}
\bibliography{main.bib}
\end{document}